\magnification=\magstep1


\def\item{\vskip1.3pt\hang\textindent}


\tolerance=300
\pretolerance=200
\hfuzz=1pt
\vfuzz=1pt

\hoffset 0cm            
\hsize=5.8 true in
\vsize=9.5 true in

\def\rightheadline{\hfil\smc\lastname\hfil\tenbf\folio}
\def\leftheadline{\tenbf\folio\hfil\smc\lastname\hfil}
\headline={\ifodd\pageno\rightheadline\else\leftheadline\fi}
\newdimen\dimenone
\def\checkleftspace#1#2#3#4#5{
 \dimenone=\pagetotal
 \advance\dimenone by -\pageshrink   
 \ifdim\dimenone>\pagegoal          
   \else\dimenone=\pagetotal
        \advance\dimenone by \pagestretch
        \ifdim\dimenone<\pagegoal
          \dimenone=\pagetotal
          \advance\dimenone by#1         
          \setbox0=\vbox{#2\parskip=0pt                
                       \hyphenpenalty=10000
                       \rightskip=0pt plus 5em
                       \noindent#3 \vskip#4}    
        \advance\dimenone by\ht0
        \advance\dimenone by 3\baselineskip
        \ifdim\dimenone>\pagegoal\vfill\eject\fi
          \else\eject\fi\fi}

\parindent=35pt
\mathsurround=1pt
\parskip=1pt plus .25pt minus .25pt
\normallineskiplimit=.99pt

\mathchardef\emptyset="001F 


%



\def\1{{\bf1}}\def\0{{\bf0}}

\def\({\bigl(}  \def\){\bigr)}
\def\<{\mathopen{\langle}}\def\>{\mathclose{\rangle}}

\def\Z{{\mathchoice{{\hbox{$\rm Z\hskip 0.26em\llap{\rm Z}$}}}%
{{\hbox{$\rm Z\hskip 0.26em\llap{\rm Z}$}}}%
{{\hbox{$\scriptstyle\rm Z\hskip 0.31em\llap{$\scriptstyle\rm Z$}$}}}{{%
\hbox{$\scriptscriptstyle\rm Z$\hskip0.18em\llap{$\scriptscriptstyle\rm Z$}}}}}}

\def\qed{\hfill {\hbox{[\hskip-0.05em ]}}}

\def\P{{\mathchoice{\hbox{$\rm I\hskip-0.14em P$}}%
{\hbox{$\rm I\hskip-0.14em P$}}%
{\hbox{$\scriptstyle\rm I\hskip-0.14em P$}}%
{\hbox{$\scriptscriptstyle\rm I\hskip-0.10em P$}}}}

\def\Q{\QQ\,\,}

\def\.{{\cdot}}
\def\|{\Vert}
\def\ssk{\smallskip}
\def\msk{\medskip}
\def\bsk{\bigskip}
\def\giantskip{\vskip2\bigskipamount}

\def\giantbreak{\par \ifdim\lastskip<2\bigskipamount \removelastskip
         \penalty-400 \giantskip\fi}

\def\nin{\noindent}
\def\cen{\centerline}
\def\pagebreak{\vskip 0pt plus 0.0001fil\break}
\def\linebreak{\break}

\def\epsilon{\varepsilon}

\font\ninerm=cmr9
\font\eightrm=cmr8
\font\sixrm=cmr6

\font\eightbf=cmbx8
\font\sixbf=cmbx6

\font\eighti=cmmi8
\font\sixi=cmmi6
\font\ninesy=cmsy9
\font\eightsy=cmsy8
\font\sixsy=cmsy6

\font\eightit=cmti8


\font\eightsl=cmsl8

\font\eighttt=cmtt8
\font\bfone=cmbx10 scaled\magstep1 
\font\bftwo=cmbx10 scaled\magstep2 
\font\smc=cmcsc10

\font\small=cmcsc8

\def\no #1. {\bigbreak\vskip-\parskip\noindent\bf #1. \quad\rm}

\def\Proposition #1. {\checkleftspace{0pt}{\bf}{Theorem}{0pt}{}
\bigbreak\vskip-\parskip\noindent{\bf Proposition #1.}
\quad\it}

\def\Theorem #1. {\checkleftspace{0pt}{\bf}{Theorem}{0pt}{}
\bigbreak\vskip-\parskip\noindent{\bf  Theorem #1.}
\quad\it}
\def\Corollary #1. {\checkleftspace{0pt}{\bf}{Theorem}{0pt}{}
\bigbreak\vskip-\parskip\nin{\bf Corollary #1.}
\quad\it}
\def\Lemma #1. {\checkleftspace{0pt}{\bf}{Theorem}{0pt}{}
\bigbreak\vskip-\parskip\noindent{\bf  Lemma #1.}\quad\it}

\def\Definition #1. {\checkleftspace{0pt}{\bf}{Theorem}{0pt}{}
\rm\bigbreak\vskip-\parskip\noindent{\bf Definition #1.}
\quad}

\def\Remark #1. {\checkleftspace{0pt}{\bf}{Theorem}{0pt}{}
\rm\bigbreak\vskip-\parskip\noindent{\bf Remark #1.}\quad}

\def\Exercise #1. {\checkleftspace{0pt}{\bf}{Theorem}{0pt}{}
\rm\bigbreak\vskip-\parskip\noindent{\bf Exercise #1.}
\quad}

\def\Example #1. {\checkleftspace{0pt}{\bf}{Theorem}{0pt}{}
\rm\bigbreak\vskip-\parskip\noindent{\bf Example #1.}\quad}
\def\Examples #1. {\checkleftspace{0pt}{\bf}{Theorem}{0pt}
\rm\bigbreak\vskip-\parskip\noindent{\bf Examples #1.}\quad}

\newcount\problemnumb \problemnumb=0
\def\Problem{\global\advance\problemnumb by 1\bigbreak\vskip-\parskip\noindent
{\bf Problem \the\problemnumb.}\quad\rm }

\def\Proof#1.{\rm\par\ifdim\lastskip<\bigskipamount\removelastskip\fi\smallskip
            \noindent {\bf Proof.}\quad}

\nopagenumbers

\def\author{}
\def\lastname{}
\def\thanks#1{\footnote*{\eightrm#1}}
\def\title{}

\def\lastname{}
\def\h{{\textstyle{1\over2}}}

\def\n{{\cal N}}

\def\text{\textstyle}
\def\disp{\displaystyle}

\def\and{{\rm and }}

\def\n{\cen{{\it W.G. Nowak}}}

\expandafter\edef\csname amssym.def\endcsname{%
       \catcode`\noexpand\@=\the\catcode`\@\space}
\catcode`\@=11
\def\undefine#1{\let#1\undefined}
\def\newsymbol#1#2#3#4#5{\let\next@\relax
 \ifnum#2=\@ne\let\next@\msafam@\else
 \ifnum#2=\tw@\let\next@\msbfam@\fi\fi
 \mathchardef#1="#3\next@#4#5}
\def\mathhexbox@#1#2#3{\relax
 \ifmmode\mathpalette{}{\m@th\mathchar"#1#2#3}%
 \else\leavevmode\hbox{$\m@th\mathchar"#1#2#3$}\fi}
\def\hexnumber@#1{\ifcase#1 0\or 1\or 2\or 3\or 4\or 5\or 6\or 7\or 8\or
 9\or A\or B\or C\or D\or E\or F\fi}

\font\tenmsb=msbm10
\font\sevenmsb=msbm7
\font\fivemsb=msbm5
\newfam\msbfam
\textfont\msbfam=\tenmsb
\scriptfont\msbfam=\sevenmsb
\scriptscriptfont\msbfam=\fivemsb
\edef\msbfam@{\hexnumber@\msbfam}
\def\Bbb#1{{\fam\msbfam\relax#1}}

\newsymbol\Bbbk 207C
\def\widehat#1{\setbox\z@\hbox{$\m@th#1$}%
 \ifdim\wd\z@>\tw@ em\mathaccent"0\msbfam@5B{#1}%
 \else\mathaccent"0362{#1}\fi}
\def\widetilde#1{\setbox\z@\hbox{$\m@th#1$}%
 \ifdim\wd\z@>\tw@ em\mathaccent"0\msbfam@5D{#1}%
 \else\mathaccent"0365{#1}\fi}
\font\teneufm=eufm10
\font\seveneufm=eufm7
\font\fiveeufm=eufm5
\newfam\eufmfam
\textfont\eufmfam=\teneufm
\scriptfont\eufmfam=\seveneufm
\scriptscriptfont\eufmfam=\fiveeufm
\def\frak#1{{\fam\eufmfam\relax#1}}

\catcode`@=11 

\expandafter\edef\csname amssym.def\endcsname{%
       \catcode`\noexpand\@=\the\catcode`\@\space}
\font\eightmsb=msbm8
\font\sixmsb=msbm6
\font\fivemsb=msbm5
\font\eighteufm=eufm8
\font\sixeufm=eufm6
\font\fiveeufm=eufm5
\newskip\ttglue
\def\eightpoint{\def\rm{\fam0\eightrm}%
  \textfont0=\eightrm \scriptfont0=\sixrm \scriptscriptfont0=\fiverm
  \textfont1=\eighti \scriptfont1=\sixi \scriptscriptfont1=\fivei
  \textfont2=\eightsy \scriptfont2=\sixsy \scriptscriptfont2=\fivesy
  \textfont3=\tenex \scriptfont3=\tenex \scriptscriptfont3=\tenex
\textfont\eufmfam=\eighteufm
\scriptfont\eufmfam=\sixeufm
\scriptscriptfont\eufmfam=\fiveeufm
\textfont\msbfam=\eightmsb
\scriptfont\msbfam=\sixmsb
\scriptscriptfont\msbfam=\fivemsb
  \def\it{\fam\itfam\eightit}%
  \textfont\itfam=\eightit
  \def\sl{\fam\slfam\eightsl}%
  \textfont\slfam=\eightsl
  \def\bf{\fam\bffam\eightbf}%
  \textfont\bffam=\eightbf \scriptfont\bffam=\sixbf
   \scriptscriptfont\bffam=\fivebf
  \def\tt{\fam\ttfam\eighttt}%
  \textfont\ttfam=\eighttt
  \tt \ttglue=.5em plus.25em minus.15em
  \normalbaselineskip=9pt
  \def\MF{{\manual opqr}\-{\manual stuq}}%
  \let\big=\eightbig
  \setbox\strutbox=\hbox{\vrule height7pt depth2pt width\z@}%
  \normalbaselines\rm}
\def\eightbig#1{{\hbox{$\textfont0=\ninerm\textfont2=\ninesy
  \left#1\vbox to6.5pt{}\right.\n@space$}}}


\csname amssym.def\endcsname


\def\al{\alpha}

\def\om{\omega}
\def\({\left(}
\def\){\right)}

\def\eq{\eqalign}

\def\klein{\eightpoint \def\smc{\small} \baselineskip=9pt}

\font\boldmas=msbm10                  
\def\Bbb#1{\hbox{\boldmas #1}}        
\def\Z{{\Bbb Z}}                        
\def\Q{{\Bbb Q}}


\font\eightrm=cmr8
\long\def\fussnote#1#2{{\baselineskip=9pt
\setbox\strutbox=\hbox{\vrule height 7pt depth 2pt width 0pt}%
\eightrm
\footnote{#1}{#2}}}
\font\boldmasi=msbm10 scaled 700      
\def\Bbbi#1{\hbox{\boldmasi #1}}      
\font\boldmas=msbm10                  
\def\Bbb#1{\hbox{\boldmas #1}}        
\def\Zi{{\Bbbi Z}}                      
\def\Pi{{\Bbbi P}}                      
\def\Qi{{\Bbbi Q}}                      



\def\dint #1 {
\quad  \setbox0=\hbox{$\disp\int\!\!\!\int$}
  \setbox1=\hbox{$\!\!\!_{#1}$}
  \vtop{\hsize=\wd1\centerline{\copy0}\copy1} \quad}

\def\drint #1 {
\qquad  \setbox0=\hbox{$\disp\int\!\!\!\int\!\!\!\int$}
  \setbox1=\hbox{$\!\!\!_{#1}$}
  \vtop{\hsize=\wd1\centerline{\copy0}\copy1}\qquad}

\def\frac#1#2{{#1\over #2}}

\def\date{\the\day.~\the\month.~\the\year}

\def\mod{\,{\rm mod}\,}
\def\klein{\eightpoint \def\smc{\small} }

\def\frac#1#2{{#1\over#2}}

\hsize=16true cm     \vsize=23.5true cm

\parindent=0cm

\def\ok{{\cal O}_K}
\def\b{{\bf b}_K}
\def\oOm{\overline{\Omega}}
\def\o#1{\overline{#1}}
\def\c{\overline{\frak c}}
\def\p{{\frak P}}
\def\n{{\cal N}}
\def\D{{\cal D}}
\def\PRO{\prod_{m,k=1\atop m\ne k}^M (a_m b_k-a_k b_m)}
\def\qed{\hfill\hbox{\vrule\vbox to 0.8em{\hsize=0.8em\hrule\noindent\
\vfill\hrule}\vrule}}
\def\Pu{\Pi_{\cup}}

\vbox{\vskip 1.2true cm}

\cen{\bftwo On the distribution of M-tuples of B-numbers}\bsk\msk

\cen{\bfone Werner Georg Nowak (Vienna)}

\vbox{\vskip 1.2true cm}

\footnote{}{\klein{\it Mathematics Subject Classification }
(2000): 11P05, 11N35.\par { \it Keywords$\,$}: 
$B$-numbers; Selberg sieve; norms of ideals in number fields.}

{\klein{\bf Abstract. } In the classical sense, the set $B$
consists of all integers which can be written as a sum of two
perfect squares. In other words, these are the values attained by
norms of integral ideals over the Gaussian field $\Qi(i)$. G.J.
Rieger (1965) and T. Cochrane / R.E. Dressler (1987) established
bounds for the number of pairs $(n,n+h)$, resp., triples
$(n,n+1,n+2)$ of $B$-numbers up to a large real parameter $x$. The
present article generalizes these investigations into two
directions: The result obtained deals with arbitrary $M$-tuples of
arithmetic progressions of positive integers, excluding the
trivial case that one of them is a constant multiple of one of the others. 
Furthermore, the estimate applies to the case of an arbitrary
normal extension $K$ of the rational field instead of $\Qi(i)$.}

\vbox{\vskip 1.2true cm}

{\bf1.~Introduction. } Already E.~Landau's in his classic
monograph [4] provided a proof of the result that the set $B$ of
all positive integers which can be written as a sum of two squares
of integers is distributed fairly regularly: It satisfies the
asymptotic formula
$$ \sum_{1\le n\le x,\ n\in B} 1 \ \sim \ {c\,x\over\sqrt{\log x}}\qquad(c>0)\,.
\eqno(1.1)$$ Almost six decades later, G.J.~Rieger [9] was the
first to deal with the question of "$B$-twins": How frequently
does it happen that both $n$ and $n+1$ belong to the set $B$? A
bit more general, he was able to show that, for any positive
integer $h$ and large real $x$,
$$ \sum_{1\le n\le x\atop n\in B,\ n+h\in B} 1 \quad \ll \quad 
\prod_{p\mid h \atop p\equiv3\mod4}\(1+{1\over p}\)\ {x\over\log x}\,. 
\eqno(1.2) $$
Later on, C.~Hooley [2] and K.-H.~Indlekofer [3], independently
and at about the same time, showed that this bound is  essentially
best possible. \ssk In 1987, T.~Cochrane and R.E.~Dressler [1]
extended the question to triples of $B$-numbers. Replacing
Rieger's sieve technique by a more recent variant of Selberg's
method, they succeeded in proving that
$$  \sum_{1\le n\le x\atop n\in B,\, n+1\in B,\, n+2\in B} 1 \quad \ll \quad
 {x\over(\log x)^{3/2}}\,. \eqno(1.3) $$ \bsk

{\bf2.~Statement of result. } In this article we intend to
generalize these estimates in two different directions: Firstly,
instead of pairs or triples we consider $M$-tuples of arithmetic
progressions $(a_m\,n+b_m)$, $m=1,\dots,M\ge2$, where
$a_m\in\Z^+$, $b_m\in\Z$ throughout. Secondly, we deal with an
arbitrary number field $K$ which is supposed to be a normal
extension of the rationals of degree $[K:\Q]=N\ge2$. Denoting by
$\ok$ the ring of algebraic integers in $K$, we put
$$ \b(n) := \cases{1& if there exists an integral ideal ${\frak A}$
in $\ok$ of norm ${\cal N}({\frak A})=n$,\cr 0 & else.} $$ Our
target is then the estimation of the sum
$$ S(x) = S(a_1,b_1, \dots ,a_M,b_M;\,x) := \sum_{1\le n\le
x}\prod_{m=1}^M \b(a_m\,n+b_m)\,. \eqno(2.1) $$ Of course, the
classic case reported in section 1 is contained in this, by the
special choice $K=\Q(i)$, the Gaussian field. \bsk

{\bf Theorem. } { \it Suppose that $(a_m,b_m)\in\Z^+\times\Z$ for
$m=1,\dots,M$, and, furthermore, $$ \prod_{m,k=1\atop m\ne k}^M
(a_m b_k-a_k b_m)\ \ne \ 0 \,. $$ Then, for large real $x$,
$$  S(a_1,b_1, \dots ,a_M,b_M;\,x) \ll \gamma(a_1,b_1, \dots
,a_M,b_M) \,{x\over(\log x)^{M(1-1/N)}}\,, $$ with
$$  \gamma(a_1,b_1, \dots ,a_M,b_M) = \prod_{p\in\Pi'}\(1+{M\over
p}\)\,, $$ the finite set of primes $\P'=\P'(a_1,b_1, \dots
,a_M,b_M)$ to be defined below in $(4.6)$. The $\ll$-constant
depends on $M$ and the field $K$, but not on $a_1,b_1, \dots
,a_M,b_M$.} \bsk\msk

{\bf 3.~Some auxiliary results.} \msk

{\bf Notation. } Variables of summation automatically range over
all integers satisfying the conditions indicated.  $p$ denotes
rational primes throughout, and $\P$ is the set of all rational
primes. $\p$ stands for prime ideals in $\ok$. For any subset
$\P^\circ\subseteq\P$, we denote by $\D(\P^\circ)$  the set of all
positive integers whose prime divisors all belong to $\P^\circ$.
The constants implied in the symbols $O(\cdot)$, $\ll$, $\gg$,
etc., may depend throughout on the field $K$ and on $M$, but not
on $a_1,b_1, \dots ,a_M,b_M$. \bsk

{\bf Lemma 1. } {\it For each prime power $p^\al$, $\al\ge1$, let
$\oOm(p^\al)$ be a set of distinct residue classes $\c$ modulo
$p^\al$. Define further
$$ \Omega(p^\al) = \{n\in\Z^+:\ n\in\bigcup_{\c\in\oOm(p^\al)}\c\ \}\,, $$
and let $$  \theta(p^\al) := 1 - \sum_{j=1}^\al {\#\oOm(p^j)\over
p^j}>0\,,\qquad \theta(1):=1\,. $$ Suppose that
$\Omega(p^\al)\cap\Omega(p^{\al'})=\emptyset$ for all primes $p$
and positive integers $\al\ne\al'$. For real $x>0$, let finally
$$  A(x) = \{n\in\Z^+:\ n\le x\quad \and\quad
n\notin\bigcup_{p\in\Pi, \al\in\Zi^+} \Omega(p^\al)\ \}\,.  $$
Then, for arbitrary real $Y>1$,
$$  \# A(x) \le {x+Y^2\over V_Y}\,, $$
where
$$  V_Y := \sum_{0<d<Y}\quad \prod_{p^\al\,\Vert\,d}
\({1\over\theta(p^\al)}-{1\over\theta(p^{\al-1})}\)\,.  $$ } \bsk

{\bf Proof. } This is a deep sieve theorem due to A.~Selberg [10].
It can be found in Y.~Motohashi [5], p.~11, and also in
T.~Cochrane and R.E.~Dressler [1]. \bsk\msk

{\bf Lemma 2. } {\it Let $(c_n)_{n\in\Zi^+}$ be a sequence of
nonnegative reals, and suppose that the Dirichlet series
$$  f(s) = \sum_{n=1}^\infty c_n\,n^{-s}  $$
converges for $\Re(s)>1$. Assume further that, for some real
constants $A$ and $\beta>0$,
$$ f(s) = (A+o(1))(s-1)^{-\beta}\,, $$
as $s\to1+$. Then, for $x\to\infty$,
$$  \sum_{1\le n\le x} {c_n\over n} \ =\
\({A\over\Gamma(1+\beta)}+o(1)\)(\log x)^{\beta}\,.  $$   } \bsk
{\bf Proof. } This is a standard Tauberian theorem. For the
present formulation, cf.~Cochrane and Dressler [1], Lemma B.
\bsk\msk

{\bf4.~Proof of the Theorem. } We recall the decomposition laws in
a normal extension $K$ over $\Q$ of degree $N\ge2$
(cf.~W.~Narkiewicz [6], Theorem 7.10.): Every rational prime $p$
which does not divide the field discriminant disc$(K)$ belongs to
one of the classes $$ \P_r = \{ p\in\P:\ (p)=\p_1 \cdots
\p_{N/r},\ {\cal N}(\p_1)=\dots=\ {\cal N}(\p_{N/r})= p^r\ \}\,,
$$ where $r$ ranges over the divisors of $N$, and
$\p_1,\dots,\p_{N/r}$ are distinct. As an easy consequence, if
$p\in\P_r$, $\al\in\Z^+$,  $$ \b(p^\al)= \cases{1 & if
$r\mid\al$,\cr 0 & else.\cr} \eqno(4.1)  $$ In order to apply
Lemma 1, we need a bit of preparation. Let
$$ \eq{\P_r^* &= \bigl\{ p\in\P_r:\ p\not| \ \prod_{m=1}^M a_m \PRO\,,
\ p\ne M-1\, \bigr\}\,, \cr \P^* &= \bigcup_{r\mid N,\ r>1}
\P_r^*\,.\cr } $$ Then we choose
$$ \oOm(p^\al) := \bigcup_{m=1}^M \left\{ \o{a_m}\,^{(-1)}
\o{(jp^{\al-1}-b_m)} :\ j=1,\dots,p-1\ \right\}\,, $$ if
$p\in\P_r^*$ and $r\not|(\al-1)$, while $\oOm(p^\al) :=\emptyset$
in all other cases. Here $\o{\,\cdot\,}$ denotes residue classes
modulo $p^\al$, in particular $\o{a_m}\,^{(-1)}$ is the class
which satisfies $\o{a_m}\,\,\o{a_m}\,^{(-1)} = \o{\,1\,}$ mod
$p^\al$. We summarize the relevant properties of these sets
$\oOm(p^\al)$, and of the corresponding sets $\Omega(p^\al)$ (see
Lemma 1), as follows. \bsk

{\bf Proposition. } {\it Suppose throughout that $p\in\P^*$ and
$\al\in\Z^+$.}\ssk (i) { \it If $p\in\P_r^*$, $r\not|(\al-1)$,
then $\oOm(p^\al)$ contains exactly $M(p-1)$ elements.} \ssk (ii)
{ \it If a positive integer $k$ lies in some $\Omega(p^\al)$, it
follows that there exists an $m\in\{1,\dots,M\}$ such that
$p^{\al-1}\parallel(a_m\,k+b_m)$.} \ssk (iii) { \it It is
impossible that there exist $m,n\in\{1,\dots,M\}$, $m\ne n$, and a
positive integer $k$, such that any $p\in\P^*$ divides both $a_m
k+b_m$ and $a_n k+b_n$.} \ssk (iv) { \it If $k\in\Omega(p^\al)$,
it follows that $$ p^{\al-1} \parallel \prod_{m=1}^M
(a_m\,k+b_m)\,.$$ Consequently,
$\Omega(p^\al)\cap\Omega(p^{\al'})=\emptyset$ for any positive
integers $\al\ne\al'$.}\ssk (v) { \it If $k\in\Omega(p^\al)$, then
$$ \prod_{m=1}^M \b(a_m\,k+b_m) = 0\,. $$ As a consequence, $$
S(x)\le \# A(x)\,,  $$ where $S(x)$ and $A(x)$ have been defined
in $(2.1)$ and Lemma 1, respectively.} \bsk

\vbox{{\bf Proof of the Proposition. } (i) Assume that two of
these residue classes would be equal, say, $\o{a_m}\,^{(-1)}
\o{(u\,p^{\al-1}-b_m)}$ and $\o{a_n}\,^{(-1)}
\o{(v\,p^{\al-1}-b_n)}$, where $u,v\in\{1,\dots,p-1\}$,
$m,n\in\{1,\dots,M\}$. Multiplying by $\o{a_m}\,\o{a_n}$, we could
conclude that $$ a_n(u\,p^{\al-1}-b_m) \equiv a_m
(v\,p^{\al-1}-b_n)\quad\mod\ p^\al\,, $$ or, equivalently, that $$
(a_n\,u-a_m\,v)p^{\al-1} \equiv a_n\,b_m-a_m\,b_n\quad\mod\
p^\al\,. \eqno(4.2)$$ Hence $p\mid(a_n\,b_m-a_m\,b_n)$, which is
only possible if $m=n$. This in turn simplifies (4.2) to
$$ a_m(u-v)p^{\al-1} \equiv 0\quad\mod\ p^\al\,, $$ thus also
$u=v$. \qed } \msk (ii) If $k\in\Omega(p^\al)$, there exist
$j\in\{1,\dots,p-1\}$, $m\in\{1,\dots,M\}$, and an integer $q$,
such that $$ a_m\,k = j\,p^{\al-1}-b_m + q\,p^\al\,. $$ From this
the assertion is obvious. \qed \msk (iii) Assuming the contrary,
we would infer that $p$ divides $$ (a_m\,k+b_m)b_n-(a_n\,k+b_n)b_m
= (a_m\,b_n-a_n\,b_m)k\,, $$ hence $p\mid k$, thus $p$ divides
also $b_m$ and $b_n$, which contradicts $p\in\P^*$. \qed \msk (iv)
This is immediate from (ii) and (iii). \qed \msk (v) By (ii),
$p^{\al-1}\parallel(a_m\,k+b_m)$ for some $m\in\{1,\dots,M\}$.
Recalling that $r\not|(\al-1)$ (otherwise $\Omega(p^\al)$ would be
empty), along with (4.1) and the multiplicativity of $\b(\cdot)$,
it is clear that $\b(a_m k+b_m)=0$. The last inequality is obvious
from the relevant definitions. \qed \bsk

We are now ready to apply Lemma 1. Choosing $Y=\sqrt{x}$ and
appealing to part (v) of the Proposition, we see that
$$  S(x) \le {2x\over V_Y}\,. \eqno(4.3) $$ To derive a lower
bound for $V_Y$, observe that $\oOm(p)=\emptyset$ for every prime
$p$, and $\#\oOm(p^2)=M(p-1)$ for each $p\in\P^*$. Further, if
$p\notin\P^*$, then $\oOm(p^j)=\emptyset$ throughout. Therefore,
if $p\in\P^*$, \ $\theta(p)=1$ and
$$ \theta(p^2)=1-{1\over p^2}\,\#\oOm(p^2) = 1- {M(p-1)\over
p^2}\,, $$ hence $$ {1\over\theta(p^2)}-{1\over\theta(p)} =
{M(p-1)\over p^2-M(p-1)} \ge {M\over p}\,.  \eqno(4.4) $$
Furthermore, for any $\al>2$,
$$ \theta(p^\al)\ge 1 - M(p-1)\sum_{2\le j\le\al}p^{-j} > 1-{M\over
p}\ge0\,,  $$ since $M(p-1)\le p^2-1$ according to clause (i) of
the Proposition, and $M=p+1$ is impossible for $p\in\P^*$. Thus
actually $\theta(p^\al)>0$ for all primes $p$ and all
$\al\in\Z^+$. Thus all the terms in the sum $V_Y$ are nonnegative,
and restricting the summation to the set
$$ Q := \{ d=d_1^2:\ d_1\in\Z^+\,,\ \mu(d_1)\ne0\,,\
d_1\in\D(\P^*)\ \}\,, $$ we conclude by (4.4) that $$ \eq{V_Y &\ge
\sum_{0<d<Y,\,d\in Q}\ \prod_{p\mid
d}\({1\over\theta(p^2)}-{1\over\theta(p)} \)\ge
\sum_{0<d_1<\sqrt{Y},\,d_1\in\D(\Pi^*)} \mu^2(d_1) \prod_{p\mid
d_1} {M\over p}=\cr &=  \sum_{0<d_1<\sqrt{Y},\,d_1\in\D(\Pi^*)}
\mu^2(d_1)\,{M^{\om(d_1)}\over d_1}\,, \cr } \eqno(4.5)  $$ where
$\om(d_1)$ denotes the number of primes dividing $d_1$. Our next
step is to take care of the primes excluded in the construction of
$\P^*$. We define $$ \P' = \P'(a_1,b_1, \dots ,a_M,b_M) :=
\bigl\{p\in\bigcup_{r\mid N\atop r>1}\P_r\, :\ p \mid \
\prod_{m=1}^M a_m \PRO\ \bigr\} \eqno(4.6) $$ and $$ \gamma =
\gamma(a_1,b_1, \dots ,a_M,b_M) := \prod_{p\in\Pi'} \(1+{M\over
p}\) = \sum_{k_1\in\D(\Pi')}\mu^2(k_1)\,{M^{\om(k_1)}\over k_1}\,.
\eqno(4.7)$$ Putting finally $$ \P_{\cup} := \bigcup_{r\mid N,\,
r>1} \P_r\,, $$ we readily infer from (4.5) and (4.7) that $$
\gamma\,V_Y\gg\sum_{k<\sqrt{Y},\,k\in\D(\Pi_{\cup})}\mu^2(k)\,
{M^{\om(k)}\over k}\,. \eqno(4.8)$$ We shall estimate this latter
sum by the corresponding generating function
$$ f(s) := \sum_{k\in\D(\Pi_{\cup})}\mu^2(k)\,{M^{\om(k)}\over k^s}
= \prod_{p\in\Pu} \(1+{M\over p^s}\) \qquad(\Re(s)>1)\,, $$
applying Lemma 2. By $h_1(s), h_2(s), \dots $ we will denote
functions which are holomorphic and bounded, both from above and
away from zero, in every half-plane $\Re(s)\ge\sigma_0>\h$. We
first observe that
$$ f(s) = h_1(s) \prod_{p\in\Pu}
\(1-p^{-s}\)^{-M}\qquad(\Re(s)>1)\,. \eqno(4.9)  $$ This follows
by a standard argument which can be found exposed neatly in
G.~Tenenbaum [11], p.~200 f. The next step is to consider the Euler
product of the Dedekind zeta-function $\zeta_K(s)$: For
$\Re(s)>1$,
$$ \eq{\zeta_K(s) &= \prod_{\p}\(1-\n(\p)^{-s}\)^{-1} = h_2(s) \prod_{r\mid
N}\(\prod_{p\in\Pi_r}\(1-p^{-rs}\)^{-N/r}\) = \cr  &=
h_3(s)\prod_{p\in\Pi_1}\(1-p^{-s}\)^{-N}\,. \cr  }   $$ Therefore,
$$ {\(\zeta(s)\)^{M}\over\(\zeta_K(s)\)^{M/N}} =  h_4(s) \prod_{p\in\Pu}
\(1-p^{-s}\)^{-M}\qquad(\Re(s)>1)\,. $$ Comparing this with (4.9),
we arrive at $$ f(s) = h_5(s)\,
{\(\zeta(s)\)^{M}\over\(\zeta_K(s)\)^{M/N}} \,. $$ From this it is
evident that, as $s\to1+$,
$$ f(s) \sim h_5(1)\,\rho_K^{-M/N}\,(s-1)^{-M+M/N}\,, $$ where
$\rho_K$ denotes the residue of $\zeta_K(s)$ at $s=1$. Lemma 2 now
immediately implies that
$$ \sum_{k<\sqrt{Y},\,k\in\D(\Pi_{\cup})}\mu^2(k)\,
{M^{\om(k)}\over k} \gg (\log Y)^{M-M/N} \gg (\log x)^{M-M/N}\,,
$$ in view of our earlier choice $Y=\sqrt{x}$. Combing this with
(4.3) and (4.8), we complete the proof of our Theorem. \qed

\bsk\msk

{\bf5.~Concluding remarks. } 1.~Taking more care and imposing
special conditions on the numbers $a_1, b_1, \dots , a_M, b_M$,
one could improve slightly on the factor $\gamma $ in our
estimate. (Observe that Rieger's bound (1.2) is in fact a bit
sharper than our general result.) But it is easy to see that
$\gamma $ is rather small anyway: By elementary facts about the
Euler totient function (see K.~Prachar [8], p.~24-28),

\vbox{$$ \gamma(a_1,b_1, \dots ,a_M,b_M)  \ll
\prod_{p\in\Pi'}\(1-{1\over p}\)^{-M} \ll (\log\log x)^M\,, $$
under the very mild restriction that, for some constant $c>0$, $$
\max_{m=1,\dots,M}\(a_m, |b_m|\) \ll \exp((\log x)^c)\,. $$} \ssk
2.~As far as the asymptotics (1.1) is concerned, the
generalization to an arbitrary normal extension $K$ of $\Q$ can be
found in W.~Narkiewicz' monograph [6], p.~361, Prop.~7.11, where
it is attributed to E.~Wirsing. For this question, the case of
non-normal extensions $K$ has been dealt with by R.W.K.~Odoni [7].
It may be interesting to extend our present problem to the
non-normal case as well. We might return to this at a later
occasion.

\vbox{\vskip 1.5true cm}

\klein

\cen{\bf References}\bsk \def\it{}

[1] T.~Cochrane \and R.E.~Dressler, {\it Consecutive triples of
sums of two squares,} Arch.~Math.~(Basel) 49 (1987), 301-304. \ssk

[2] C.~Hooley, On the intervals between numbers that are sums of
two squares, III, J.~Reine Angew.~Math. 267 (1974), 207-218. \ssk

[3] K.-H.~Indlekofer, Scharfe untere Absch\"atzung f\"ur die
Anzahlfunktion der $B$-Zwillinge, Acta Arith. 26 (1974), 207-212.
\ssk

[4] E.~Landau, Handbuch der Lehre von der Verteilung der
Primzahlen, Leipzig, Teubner, 1909.  \ssk

[5] Y.~Motohashi, {\it Lectures on sieve methods and prime number
theory,} Tata Institute of Fund.~Research, Bombay, 1983.  \ssk

[6] W.~Narkiewicz, Elementary and analytic theory of algebraic
numbers, 2nd ed., Berlin 1990. \ssk

[7] R.W.K.~Odoni, On the norms of algebraic integers, Mathematika
22 (1975), 71-80. \ssk

[8] K.~Prachar, Primzahlverteilung, Berlin-G\"ottingen-Heidelberg,
1957. \ssk

[9] G.J.~Rieger, {\it Aufeinanderfolgende Zahlen als Summen von
zwei Quadraten,} Indag.~Math. 27 (1965), 208-220.  \ssk

[10] A.~Selberg,  {\it Remarks on multiplicative functions,}
Number Theory Day, Proc.~Conf., New York 1976, Lecture Notes in
Math. 626 (1977), 232-241.  \ssk

[11] G.~Tenenbaum, Introduction to analytic and probabilistic
number theory, Cambridge 1995.

\vbox{\vskip2true cm}

\vbox{Werner Georg Nowak \par Institute of Mathematics \par
Department of Integrative Biology \par Universit\"at f\"ur
Bodenkultur Wien \par Gregor Mendel-Stra\ss e 33 \par 1180 Wien,
Austria \ssk E-mail: {\tt nowak@mail.boku.ac.at}\ssk Web:
http://www.boku.ac.at/math/nth.html}

\bye